\documentclass[psamsfonts]{pspum-l}
\usepackage{amsmath, amsfonts, stmaryrd, eucal}

\newtheorem{theorem}{Theorem}[section]

\theoremstyle{definition}

\theoremstyle{remark}

\numberwithin{equation}{section}

\newcommand{\calB}{\mathcal{B}}
\newcommand{\calD}{\mathcal{D}}
\newcommand{\calL}{\mathcal{L}}
\newcommand{\calR}{\mathcal{R}}
\newcommand{\calO}{\mathcal{O}}
\newcommand{\frakX}{\mathfrak{X}}
\newcommand{\gotho}{\mathfrak{o}}
\newcommand{\AAA}{\mathbb{A}}
\newcommand{\CC}{\mathbb{C}}
\newcommand{\PP}{\mathbb{P}}
\newcommand{\QQ}{\mathbb{Q}}
\newcommand{\ZZ}{\mathbb{Z}}
\newcommand{\del}{\partial}
\newcommand{\bB}{\mathbf{B}}
\newcommand{\bv}{\mathbf{v}}

\DeclareMathOperator{\coh}{coh}
\DeclareMathOperator{\crys}{crys}
\DeclareMathOperator{\dR}{dR}
\DeclareMathOperator{\et}{et}
\DeclareMathOperator{\Fil}{Fil}
\DeclareMathOperator{\Frac}{Frac}

\DeclareMathOperator{\Gal}{Gal}
\DeclareMathOperator{\GL}{GL}
\DeclareMathOperator{\HT}{HT}
\DeclareMathOperator{\inte}{inte}

\DeclareMathOperator{\rig}{rig}
\DeclareMathOperator{\Spec}{Spec}
\DeclareMathOperator{\speci}{sp}
\DeclareMathOperator{\Spf}{Spf}
\DeclareMathOperator{\st}{st}

\begin{document}

\title{$p$-adic Cohomology}

\author{Kiran S. Kedlaya}
\address{Department of Mathematics, Massachusetts Institute of Technology,
77 Massachusetts Avenue, Cambridge, MA 02139}
\email{kedlaya@math.mit.edu}
\urladdr{http://math.mit.edu/\~{}kedlaya/}
\thanks{The author was supported in part by NSF Grant DMS-0400747.}

\date{April 16, 2008.}


\maketitle

The purpose of this paper is to survey some recent results in the theory
of ``$p$-adic cohomology'', by which we will mean several different
(but related) things:
the de Rham or $p$-adic \'etale cohomology of varieties over $p$-adic fields,
or the rigid cohomology of varieties over fields of characteristic $p>0$.
Our goal is to update Illusie's beautiful 1994 survey
\cite{illusie} by reporting on some of the many interesting results 
that postdate it. In particular, we concentrate more on the present
and near-present than the past; \cite{illusie} provides a much better
historical background than we can aspire to, and is a more advisable starting
point for newcomers.

Before beginning, it is worth pointing out one (but not the only) crucial
reason why so much progress has been made since the appearance of
\cite{illusie}. In the mid-1990s, it suddenly became possible to circumvent
the resolution of singularities problem in positive characteristic thanks to
de Jong's alterations theorems \cite{dejong}, which provide forms of
``weak resolution'' and 
``weak semistable reduction''.
(The weakness is the introduction of an unwanted but often relatively 
harmless finite extension of the function field: for instance, an
alteration is a morphism which is proper, dominant, and generically finite
rather than generically an isomorphism.)
The importance of de Jong's results cannot be overstated;
they underlie almost every geometric argument cited in this paper! See
\cite{berthelot-alter} for more context regarding alterations.

\section{Rigid cohomology}

For concepts in rigid analytic geometry which we do not explain here, we recommend
\cite{fresnel-vanderput} (or for more foundational details \cite{bgr}).

\subsection{Preview: Monsky-Washnitzer cohomology}

Let $k$ be a field of characteristic $p>0$, and let $K$ be a complete
discretely valued field of characteristic $0$ with residue field $k$
and ring of integers $\gotho$,
equipped with a lift of some power of the absolute Frobenius on $k$.
For instance, one may take $k$ perfect and $K = \Frac W(k)$, in which case
the lift of Frobenius exists and is unique.

Monsky and Washnitzer \cite{monsky-washnitzer, monsky1, monsky2}
(see also \cite{vanderput})
described a cohomology theory for smooth affine $k$-varieties,\footnote{For 
simplicity, here ``$k$-variety'' will mean ``reduced separated scheme
of finite type over $k$''; some results cited do not actually depend on
the reduced or separated
hypotheses.} which they called ``formal cohomology''.
Suppose $A = \gotho[x_1,\dots,x_n]/I$ is a smooth affine $\gotho$-algebra.
Define the ring of \emph{overconvergent power series}
$\gotho \langle x_1, \dots, x_n \rangle^\dagger$ as the union, over
all $\rho>1$, of the ring of power series in 
$\gotho \llbracket x_1, \dots, x_n \rrbracket$ which converge
for $|x_1|, \dots, |x_n| \leq \rho$; put 
\[
A^\dagger = \gotho \langle x_1, \dots, x_n \rangle^\dagger / 
I \gotho_K \langle x_1, \dots, x_n \rangle^\dagger.
\]
The Monsky-Washnitzer cohomology of $\Spec (A \otimes_{\gotho} k)$
is then the cohomology of the complex
of continuous differentials of $A^\dagger \otimes_{\gotho} K$.

The bad news about this construction is that 
the above description leaves many questions unanswered: whether
a smooth $k$-variety can be lifted to $\gotho$ (so that the above construction
is possible), whether it is independent of the choice of lift, whether
maps between smooth $k$-varieties can be lifted, and whether the induced
maps on cohomology are functorial (and in particular independent of
the choice of lift). All of these have affirmative answers, but with a bit
of work required; see \cite{vanderput}.

The good news is how simple the construction is to describe. This
makes it the centerpiece
of much of the theoretical analysis of rigid cohomology, but it also has
a quite unexpected side benefit: the construction has attracted much interest
from computational applications
that require the determination of the zeta function or related 
information\footnote{The ``related information'' is often the order of a
Jacobian group, but sometimes not. One example:
Mazur-Stein-Tate \cite{mst} use the matrix of the Frobenius action on
$p$-adic cohomology to compute
$p$-adic global canonical heights 
on elliptic curves over $\QQ$. Another example: Voloch
and Zarzar \cite{voloch-zarzar} use upper bounds on Picard numbers
of surfaces to construct good error-correcting codes.}
from a given variety over a finite field. For instance, this occurs in
cryptography based on elliptic or hyperelliptic curves.
The idea to use $p$-adic cohomological methods for such computations is due
to Wan \cite{lauder-wan, wan-msri}; see \cite{kedlaya-compzeta} for a survey
of some subsequent developments.

\subsection{Construction of rigid cohomology}

Monsky-Washnitzer cohomology can be thought of as an analogue of algebraic
de Rham cohomology for smooth affine varieties. Berthelot 
\cite{berthelot-memoires} (also see \cite{berthelot-finitude})
realized that
this should be generalized to an analogue of the algebraic de Rham cohomology
of arbitrary varieties, constructed (following Herrera-Lieberman 
\cite{herrera-lieberman} and Hartshorne \cite{hartshorne1, hartshorne2})
by locally embedding a given variety into a smooth $k$-variety which lifts
to $\gotho_K$.

To be specific, suppose that $X$ is a $k$-variety, $X \hookrightarrow Y$
is an open immersion of $k$-varieties with $Y$ proper over $k$, 
and $Y \hookrightarrow P_k$ is a closed
immersion for $P$ a smooth formal $\gotho_K$-scheme. Let $P_K$
denote the Raynaud generic fibre of $P$; its points correspond to integral
formal
subschemes of $P$ which are finite flat over $\gotho_K$. In particular, $P_K$
admits a specialization map $\speci: P_K \to P_k$; for a subset $U$
of $P_k$, we write
$]U[$ for $\speci^{-1}(U)$ and call it the \emph{tube} of $U$ in $P_K$.

A \emph{strict neighborhood} of $]X[$ in $]Y[$ is an admissible open subset
of $]Y[$ which together with $]Y \setminus X[$ forms an admissible covering of
$]Y[$. The \emph{rigid cohomology} $H^i_{\rig}(X/K)$
of $X$ is then constructed as the direct
limit of the de Rham cohomologies of strict neighborhoods of $]X[$ in $]Y[$.
There is a related but slightly more complicated construction of
\emph{rigid cohomology with compact supports} $H^i_{c,\rig}(X/K)$.

As for Monsky-Washnitzer cohomology, one must make some laborious calculations
to verify that the construction of rigid cohomology with and without supports
are independent of choices and appropriately functorial.
We note in passing that it would simplify foundations\footnote{For instance,
either of these descriptions might make it easier to consider rigid cohomology
for algebraic stacks, which should compare to the crystalline cohomology
for algebraic stacks described in \cite{olsson-crys}. For simplicity, we
withhold any further discussion of stacks from this paper.}
to have either
a description of rigid cohomology either in terms of cohomology on an
appropriate site, or a suitable de Rham-Witt complex; these have been
developed by le Stum \cite{lestum} and Davis, Langer, and Zink
\cite{davis-langer-zink}, respectively.

Substituting de Jong's alterations theorem for resolution of singularities
in a program suggested in \cite[\S 4.3]{illusie},
and using a comparison theorem between rigid cohomology and
crystalline cohomology (the latter being described in \cite{berthelot-thesis, 
berthelot-ogus1}),
Berthelot \cite{berthelot-finitude, berthelot-poincare}
established the following results.
\begin{enumerate}
\item[(a)]
The vector spaces $H^i_{\rig}(X/K)$ for $X$ smooth, and
$H^i_{c,\rig}(X/K)$ for $X$ arbitrary, are finite dimensional
over $K$ 
\item[(b)]
For $X$ smooth of dimension $d$,
there is a perfect Poincar\'e duality 
 pairing
\[
H^i_{\rig}(X/K) \times H^{2d-i}_{c,\rig}(X/K) \to K.
\]
\item[(c)]
For $X_1, X_2$ smooth, there is a K\"unneth decomposition
\[
H^{i}_{\rig}(X_1 \times_k X_2/K) \cong 
\oplus_j H^j_{\rig}(X_1/K) \otimes_K H^{i-j}_{\rig}(X_2/K);
\]
for $X_1, X_2$ arbitrary, the analogous decomposition holds for
cohomology with compact supports.
\end{enumerate}
Finite dimensionality for general $X$ was deduced (using Berthelot's
work) by Grosse-Kl\"onne
\cite{grosse-klonne2} using ``dagger spaces'' as introduced in
\cite{grosse-klonne}; these are a rigid analogue of Meredith's weak formal
schemes \cite{meredith}.

In addition, the existence and basic properties of cycle class
maps have been established by Petrequin \cite{petrequin}. This means
that rigid cohomology is indeed a Weil cohomology in the sense of 
Kleiman \cite{kleiman}.

\subsection{Overconvergent $F$-isocrystals} 
\label{subsec:overcon}

A natural next step after establishing the basic properties of rigid
cohomology is to look for an appropriate category of coefficient objects.
One natural category are the \emph{convergent $F$-isocrystals};\footnote{This
somewhat awkward name, and its ``overconvergent'' sibling, deserve some
clarification. The ``convergence'' here is of the Taylor series isomorphism;
the $F$ denotes the Frobenius action; the ``isocrystal'' is short for
``crystal up to isogeny'', which is how these objects first arose in
the work of Berthelot and Ogus \cite{berthelot-ogus2, ogus1, ogus2}.} in the
notation of the previous section, these are coherent modules with connection
on the tube $]X[$ which induce Taylor isomorphisms between the two pullbacks
to the tube of $X$ in $(P \times_{\gotho_K} P)_K$, equipped
with Frobenius action.

To obtain better cohomological properties,\footnote{One way to envision the
difference between
convergence and overconvergence from a geometric viewpoint, following
a suggestion of Daqing Wan, is that
convergent and overconvergent $F$-isocrystals correspond to motives
with coefficients in finite extensions of
$\ZZ_p$ and $\ZZ$, respectively.} we must restrict attention to the 
\emph{overconvergent $F$-isocrystals}.
These are coherent modules with connection
on a strict neighborhood of the tube $]X[$, 
which induce isomorphisms between the two pullbacks
to a strict neighborhood of the tube of $X$ in $(P \times_{\gotho_K} P)_K$,
again equipped with Frobenius action.
There is a natural faithful functor from overconvergent $F$-isocrystals
to convergent $F$-isocrystals; on smooth varieties,\footnote{The same result
for general varieties can probably be deduced using a descent argument, but
this does not seem to have been verified.} this functor is fully
faithful \cite{kedlaya-full}.
The rigid cohomology of $X$ with coefficients in an overconvergent $F$-isocrystal
can be defined, in the local situation, as the direct limit of the
de Rham cohomology of the connection module over strict neighborhoods of $]X[$;
again, there is a similar definition of cohomology with compact supports.

The rigid cohomology of overconvergent $F$-isocrystals on curves was
closely analyzed by Crew \cite{crew-finiteness}, who proved finite
dimensionality and Poincar\'e duality, under a certain hypothesis
which has since been verified; see Section~\ref{subsec:localmono}.
One interesting aspect of Crew's work is its blend of ideas from
algebraic geometry and functional analysis, which played a crucial
role in subsequent developments.

Based on Crew's work, Kedlaya\footnote{This argument is an exception
to our initial comment about de Jong's alterations theorem; here one
also uses a geometric argument \cite{kedlaya-etale}
based on a higher-dimensional analogue of
``Abhyankar's trick'', a strong form of Belyi's theorem in positive
characteristic.}
 \cite{kedlaya-finiteness} established
analogues of Berthelot's finiteness, Poincar\'e duality, and K\"unneth
formula results for overconvergent $F$-isocrystals.
Again it takes more work to obtain finiteness of cohomology without
supports on nonsmooth schemes; for this one needs cohomological
descent for proper hypercoverings, developed by Chiarellotto-Tsuzuki
\cite{chiarellotto-tsuzuki, tsuzuki1, tsuzuki2}.

Some other results
 have been successfully
analogized into rigid cohomology.
The Grothendieck-Ogg-Shafarevich formula, expressing Euler characteristics
of overconvergent $F$-isocrystals on
curves in terms of certain Swan conductors, follows from a local
index theorem of
Christol-Mebkhout \cite{christol-mebkhout4} and a ``Swan conductor
equals irregularity''\footnote{The analogy between irregularity and 
Swan conductors has also been pointed out in the setting of
complex analytic de Rham cohomology, by Bloch and Esnault
\cite{bloch-esnault}.} theorem of Crew \cite{crew-canonical},
Matsuda \cite{matsuda}, and Tsuzuki \cite{tsuzuki-swan}.
The Lefschetz formula for Frobenius in rigid
cohomology was established by \'Etesse and le Stum \cite{etesse-lestum}.
An analogue of Deligne's ``Weil II'' purity theorem, using a version of
Laumon's Fourier transform (see Section~\ref{subsec:d-mod}), was established by
Kedlaya \cite{kedlaya-weilii}, building on work of Crew
\cite{crew-mono, crew-finiteness}; an extension to certain complexes
of arithmetic $\mathcal{D}$-modules (see Section~\ref{subsec:d-mod})
was given by Caro \cite{caro-devissage}.

One can ask a seemingly limitless number
of further questions asking for an analogue of a given result
in \'etale cohomology; most of these remain as yet unconsidered.
One particularly intriguing one is Mokrane's analogue
\cite{mokrane} of the weight-monodromy conjecture in \'etale cohomology;
it might be possible to interpret this question in terms of rigid
geometry. If one can relate the problem to Grosse-Kl\"onne's dagger spaces,
at least in the setting of semistable reduction
(which is specially treated in \cite{grosse-klonne2}), one might
obtain results stronger than, not just equal to, those
known in \'etale cohomology. Such a geometric interpretation has already
been given for curves by Coleman and Iovita \cite{coleman-iovita-old,
coleman-iovita}.

\subsection{Local monodromy in rigid cohomology}
\label{subsec:localmono}

In this section, we expand on the notion of ``local monodromy'' in rigid
cohomology, which underpins most of the results on overconvergent 
$F$-isocrystals
cited in the previous section. It also exhibits a surprising link
with $p$-adic Hodge theory; see Section~\ref{subsec:analytic}.
A nice discussion of this topic is given by Colmez 
\cite{colmez-monodromie}.

In Crew's analysis \cite{crew-finiteness} of the cohomology of an overconvergent
$F$-isocrystal on a curve, one is led to consider what amounts to local
cohomology at the missing points. A missing point\footnote{One
can make a similar analysis for non-rational closed points, which we omit
for simplicity.} over $k$
lifts to an open unit disc over $K$, and the $F$-isocrystal only extends
to an unspecified annulus near the boundary of that disc. One thus 
naturally obtains a finite free module over the \emph{Robba ring}
$\calR_K$ of series $\sum_{n \in \ZZ} c_n t^n$ convergent on some 
unspecified annulus $\eta < |t| < 1$, equipped with commuting Frobenius
and connection actions.

This setup led Crew to propose (modulo a later refinement by Tsuzuki
\cite{tsuzuki-slope} and a reformulation by de Jong 
\cite{dejong-icm}) the following then-conjectural statement. 
Let $\calR^{\inte}_K$ denote the subring of $\calR_K$ of series with
coefficients in $\gotho_K$; this is a noncomplete but henselian discrete
valuation ring.

\begin{theorem}[Local monodromy theorem in rigid cohomology]
 \label{thm:plmt}
Let $M$ be a finite free module over the Robba ring $\calR_K$,
equipped with commuting Frobenius and connection actions.
Then $M$ admits a filtration stable under the actions, 
whose successive quotients become trivial connection modules
after tensoring over $\calR^{\inte}_K$ with some unramified extension.
\end{theorem}

Since there are two essential structures attached to $M$, the Frobenius
and connection actions, it is fitting that there are two approaches to
proving Theorem~\ref{thm:plmt} emphasizing the two structures. One approach
uses the theory of $p$-adic differential equations initiated by Dwork and 
Robba, specifically the $p$-adic local index theory of Christol-Mebkhout
\cite{christol-mebkhout1, christol-mebkhout2, christol-mebkhout3,
christol-mebkhout4}. That theory makes essentially no reference
to Frobenius actions; it is combined with further analysis in distinct ways
by Andr\'e \cite{andre} and Mebkhout \cite{mebkhout} to give two proofs
of Theorem~\ref{thm:plmt}.

The other approach, following a strategy suggested by Tsuzuki \cite{tsuzuki-slope},
is to develop a structure theorem for Frobenius actions on the Robba ring,
related to the Dieudonn\'e-Manin classification of rational Dieudonn\'e
modules.
Such a theorem is given by Kedlaya \cite{kedlaya-local} (see also
\cite{kedlaya-slope} for a simplified presentation); it combines with 
an analysis of Theorem~\ref{thm:plmt} in the case of ``unit-root Frobenius''
due to Tsuzuki  \cite{tsuzuki-unitroot} to again yield Theorem~\ref{thm:plmt}.

It is an interesting problem to extend Theorem~\ref{thm:plmt} to a sensible
notion of local monodromy for overconvergent $F$-isocrystals on 
higher-dimensional spaces. The ``generic local monodromy theorem'' of
\cite{kedlaya-finiteness} makes a first run at this; a more complete 
answer would follow from resolution of a conjecture of Shiho \cite{shiho2}, 
which asserts that any overconvergent $F$-isocrystal can be pulled back along
a suitable alteration to obtain something extendable to a log-isocrystal
on a proper variety. For more discussion of this question,
see \cite{kedlaya-semi1} (where it is referred to as the ``semistable reduction
problem'' for overconvergent $F$-isocrystals).

\section{Arithmetic $\mathcal{D}$-modules}

See
Berthelot's excellent survey \cite{berthelot-dmod-intro}.
for a more detailed description of most of what we explain here, except
for Caro's work which postdates \cite{berthelot-dmod-intro}.

\subsection{$\mathcal{D}$-modules}
\label{subsec:d-mod}

The category of overconvergent $F$-isocrystals, while useful in many
ways, suffers the defect of not supporting Grothendieck's six operations;
that is because this category is only analogous to the category of lisse
(smooth) sheaves in $\ell$-adic cohomology. A proper category of
coefficient objects would also include analogues of constructible sheaves;
in de Rham cohomology, this is accomplished by considering the desired
coefficient objects to be modules for a typically noncommutative ring
of differential operators, as in the work of Bernstein, Kashiwara, Mebkhout,
etc.

Motivated by this last consideration (and by some initial $p$-adic 
constructions of Mebkhout and Narv\'aez-Macarro \cite{mnm}), 
Berthelot \cite{berthelot-dmod1,
berthelot-dmod2}
has introduced a notion of \emph{arithmetic $\mathcal{D}$-modules}.
In the smooth liftable case these admit a ``Monsky-Washnitzer''
description: suppose
$\frakX$ is a smooth formal scheme over $\gotho_K$, and let
$\calD_{\frakX}$ be the usual ring sheaf 
of algebraic differential operators on
$\frakX$ (with respect to $\gotho_K$, but to lighten notation we suppress this). 
For $m$ a positive integer,
if $x_1, \dots, x_n$ are local coordinates, then the ring subsheaf of
$\calD_{\frakX}$ 
generated by the operators 
$\frac{1}{(p^i)!} \frac{\del^{p^i}}{\del x_j^{p^i}}$ for $i \leq m$ and $j=1,
\dots,n$ does not depend on the choice of coordinates; the description thus
sheafifies to give a ring subsheaf $\calD_{\frakX}^{(m)}$ of
$\calD_{\frakX}$. Let $\widehat{\calD}_{\frakX}^{(m)}$
be the $p$-adic completion of $\calD_{\frakX}^{(m)}$, and put
\[
\calD_{\frakX}^\dagger = \bigcup_{m=1}^\infty 
\widehat{\calD}_{\frakX}^{(m)},
\qquad
\calD_{\frakX, \QQ}^\dagger =
\calD_{\frakX}^\dagger \otimes_{\ZZ} \QQ.
\]
Berthelot showed that $\calD_{\frakX,\QQ}^\dagger$ is a coherent sheaf of 
rings, and that its category of coherent modules is functorial in the special
fibre $X = \frakX_k$; this category includes the convergent isocrystals
on $X$.

In practice, one proves foundational
results about objects associated to $\calD_{\frakX,\QQ}^\dagger$,
like the derived category $D^b_{\coh}(\calD_{\frakX,\QQ}^\dagger)$
of bounded complexes of $\calD_{\frakX,\QQ}^\dagger$-modules
 with coherent cohomology, by writing them as direct limits of related
objects over the $\widehat{\calD}_{\frakX}^{(m)}$ and proving the
statements there by working modulo a power of $p$. For instance,
this is how one constructs the standard cohomological operations for
$\calD$-modules, namely internal and external tensor product,
direct and inverse image, and exceptional
direct and inverse image \cite[\S 4.3]{berthelot-dmod-intro}.
It is also how one establishes the $\calD$-module version of Serre
duality, in this setting due to Virrion \cite{virrion1, virrion2, virrion3,
virrion4, virrion5}; its compatibility with Frobenius is due to Caro
\cite{caro-frobrel}.

\subsection{Overconvergent singularities}

One might think that since $\calD_{\frakX,\QQ}^\dagger$ is a ring
of overconvergent differential operators, we should be able to use it
to talk about overconvergent isocrystals. That is not the case, though, because
$\calD_{\frakX,\QQ}^\dagger$ is only ``overconvergent in the
differential direction'' and not in the ``coordinate direction''
(since in particular $\calO_\frakX \subset \calD_{\frakX,\QQ}^\dagger$).
To fix this,\footnote{One could presumably
also fix this by redoing the theory with $\frakX$ taken to be a 
Meredith weak formal
scheme. Indeed, this is the original approach of Mebkhout and 
Narv\'aez-Macarro
in \cite{mnm}; however, Berthelot's approach is better suited for relating
the theory to the rigid cohomology of nonsmooth varieties.} 
one allows consideration of differential operators with
overconvergent singularities along a divisor, as follows.
(For a more detailed description, see \cite[\S 4.4]{berthelot-dmod-intro}.)

Let $Z$ be a divisor on $X = \frakX_k$. Since we are giving a local 
description,
we will assume $\frakX = \Spf A$ is affine and that $Z$ is cut out within $X$
by the reduction
of some $f \in \Gamma(\calO_{\frakX}, \frakX)$. 
Let $\widehat{\calB}_{\frakX}^{(m)}(Z)$ be the ring sheaf corresponding
to the completion of $A [ T] / (f^{p^{m+1}} T - p)$.
Change the subscript $\frakX$ to $\frakX,\QQ$ 
to denote tensoring with $\QQ$ over $\ZZ$.
Put
\[
\calO_{\frakX,\QQ}({}^\dagger Z) = \varinjlim_{m}
\widehat{\calB}_{\frakX,\QQ}^{(m)}(Z);
\]
this is the sheaf of functions on $\frakX$ with overconvergent singularities
along $Z$. An important result is that $\calO_{\frakX,\QQ}({}^\dagger Z)$
is a coherent $\calD^\dagger_{\frakX, \QQ}$-module 
\cite[Th\'eor\`eme~4.4.7]{berthelot-dmod-intro}.
To get a sheaf of overconvergent differential operators with
overconvergent singularities, we take a direct limit of completed tensor
products: namely, define
\[
\calD^\dagger_{\frakX, \QQ}({}^\dagger Z) =
\varinjlim_{m} 
(\widehat{\calB}_{\frakX}^{(m)}(Z) \widehat{\otimes}_{\calO_{\frakX}}
\widehat{\calD}_{\frakX}^{(m)}) \otimes_{\ZZ} \QQ.
\]
One can carry over the study of cohomological operations and duality
to this setting; see the Virrion and Caro references above.

For $\frakX$ proper, $\calO_{\frakX,\QQ}({}^\dagger Z)$
gives a sheafified version of the 
Monsky-Washnitzer algebra associated to $U = X \setminus Z$.
This gives a way to equip an overconvergent isocrystal on $U$
with the structure of a $\calD^\dagger_{\frakX,\QQ}({}^\dagger Z)$-module,
whose cohomology is directly 
related to the rigid cohomology of the isocrystal.
This $\calD^\dagger_{\frakX,\QQ}({}^\dagger Z)$-module is in fact
coherent, but this requires a nontrivial argument,
given by Caro \cite{caro-surcon}.

\subsection{Fourier transforms}

One important component of the study of arithmetic $\calD$-modules is
the Fourier transform; it analogizes both the natural Fourier transform
in the algebraic and analytic $\calD$-module settings (more on which shortly)
and the geometric Fourier transform introduced by Deligne and Laumon
in \'etale cohomology. In this context, the Fourier transform was constructed
by Huyghe in her thesis \cite{huyghe-thesis}, the contents of which
appear in a series of papers \cite{huyghe-acyc, huyghe-interp,
huyghe-trans, huyghe-affinite, huyghe-affinite2, huyghe-dwork}.

The Fourier transform\footnote{We are describing 
the Fourier transform without supports;
there is also a Fourier transform with compact supports, but Huyghe
\cite{huyghe-thesis} showed using duality that as in the
analogous settings, the two transforms coincide.}
 is easiest to describe on the affine line, so
let us do that now. Assume that $K$ contains a chosen root $\pi$ of the
equation $\pi^{p-1} = -p$.
Take $\frakX = \PP^1_{\gotho_K}$, $Z$ to be the
point at infinity, and let $x$ be the coordinate on $\AAA^1_{K}$. 
Let $\calL_\pi$ denote the \emph{Dwork isocrystal}
on $\AAA^1_k$: it corresponds to a 
$\calD^\dagger_{\frakX,\QQ}(\infty)$-module
which over $\calO_{\frakX,\QQ}({}^\dagger Z)$ is free of rank 1 generated
by $\bv$, such that the action of $\frac{d}{dx}$ on $\bv$ is multiplication
by $\pi$. Then the Fourier transform on $\calD^{\dagger}_{\frakX,\QQ}$-modules
is defined by appropriately interpreting
the instruction: ``Construct an integral operator with kernel $\calL_\pi$.''
Loosely, one pulls back from $\frakX$ to one factor of
$\frakX \times \frakX$, tensors
with the pullback of $\calL_{\pi}$ along the multiplication map
$\AAA^1 \times \AAA^1$, then pushes forward along the other factor.
(The reality is a bit more complicated because one must work on
$\PP^1$ and not $\AAA^1$, where there is a bit of blowing up involved.
See \cite[\S 4]{huyghe-thesis} or \cite[\S 3]{huyghe-dwork}.)

This construction happily admits an alternate description which clarifies
some of its properties. Write $\del^{[i]}$ for $\frac{1}{i!} \del^i$. Let
$A(K)$ denote\footnote{In Huyghe's notation, this ring would be denoted
$A_1(K)$, as it has an $n$-dimensional analogue $A_n(K)$.}
the set of formal sums $\sum_{i,j=0}^\infty c_{i,j} x^i \del^{[j]}$
such that the ordinary power series
$\sum_{i,j=0}^\infty c_{i,j} x^i \del^j$
belongs to $K \langle x,\del \rangle^\dagger$.
In fact, $A(K)^\dagger$ is the set of global sections of
$\calD^\dagger_{\frakX,\QQ}({}^\dagger Z)$; in particular, 
$A(K)^\dagger$ admits a
noncommutative ring structure in which 
$\del x - x \del = 1$.
Moreover, any coherent $\calD^\dagger_{\frakX,\QQ}({}^\dagger Z)$-module
is represented by a coherent $A(K)^\dagger$-module \cite{huyghe-interp}.
In this notation, the Fourier transform of a coherent $A(K)^\dagger$-module
is given, up to a shift in degree, by pullback by the map
$A(K)^\dagger \to A(K)^\dagger$ defined by $x \mapsto - \del/\pi$,
$\del \mapsto \pi x$ \cite{huyghe-trans}.

The Fourier transform for arithmetic $\calD$-modules is expected to
have the same sorts of uses as in \'etale cohomology.
At least one of these has already been realized:
Mebkhout \cite{mebkhout1}
suggested using a $p$-adic Fourier transform on the affine line to implement
Laumon's proof of Deligne's ``Weil II'' in $p$-adic cohomology,
and this has been done (see Section~\ref{subsec:overcon}).

\subsection{Holonomicity: unfinished business}

The concept of an arithmetic
$\calD$-module is something of a hybrid, sharing characteristics
both of the complex analytic and the algebraic analogues. In particular,
for a coherent $\calD_{\frakX,\QQ}^\dagger$-module equipped
with an action of absolute Frobenius (i.e., a coherent 
$F$-$\calD_{\frakX,\QQ}^\dagger$-module), one can define the
characteristic cycle, prove an analogue of Bernstein's inequality, and define
a \emph{holonomic} $F$-$\calD_{\frakX,\QQ}^\dagger$-module as one
in which equality holds in Bernstein's inequality. However, the definition of
holonomicity
cannot be made as in the algebraic case
because we are working with differential operators
of infinite order, and it cannot be made as in the analytic case because
we do not\footnote{This is not to say categorically that such an analogue 
does not 
exist! A truly analytic approach to the notion of holonomicity would be quite
interesting.} have an analogue of pseudodifferential operators and their symbols.
Instead, the notion of holonomicity is defined using a process of Frobenius
descent: the presence of a Frobenius structure makes it possible to descend
a module from $\calD^\dagger_{\frakX,\QQ}$ to $\widehat{\calD}_{\frakX}^{(m)}$
for some $m$. One can then descend further to remove the completion, 
finally ending up in an algebraic situation where one can resort to the
usual concept of a good filtration in order to define the characteristic
variety. This process is described in detail in \cite[\S 5]{berthelot-dmod-intro}.

The good news is that this notion of holonomicity leads to a natural
proposal for a category of coefficients in $p$-adic cohomology (by taking
the derived category of bounded complexes with holonomic cohomology plus
some restriction on supports).
The bad news is that the indirect nature of the definition of holonomicity
makes it extremely difficult to verify the stability of holonomicity under
the cohomological operations! As a result, Berthelot's program to construct
a good $p$-adic coefficient theory has been stalled; see 
\cite[\S 5]{berthelot-dmod-intro} for several conjectures which await resolution.

A possible route around this difficulty has been proposed by Caro
\cite{caro-surhol}, who defines
a category of \emph{overholonomic} $F$-$\calD_{\frakX,\QQ}^\dagger$-modules,
by building into the definition a certain stability of coherence
under smooth base change, and the ability to perform d\'evissages in
overconvergent $F$-isocrystals in the manner of \cite{caro-devissage}.
 Caro shows that this category is stable under all of
the cohomological operations except internal and external tensor product.
He also shows that Berthelot's conjectures imply the coincidence between
the notions of holonomicity and overholonomicity (and of a related notion
of ``overcoherence'' introduced in \cite{caro-surcoh}). Moreover,
he shows that unit-root overconvergent $F$-isocrystals are overholonomic,
as are arbitrary overconvergent $F$-isocrystals on a curve 
\cite{caro-courbes}.

It is worth pointing out that the last two assertions of the previous paragraph
are based on cases of Shiho's conjecture on logarithmic extensions of
overconvergent $F$-isocrystals (see Section~\ref{subsec:localmono}). That
conjecture is a theorem of Tsuzuki for unit-root isocrystals
\cite{tsuzuki-duke} and of Kedlaya for curves \cite{kedlaya-semicurve};
it would appear that (by adapting Caro's arguments to a suitable logarithmic
setting) a resolution of Shiho's conjecture could be used to 
establish 
overholonomicity of arbitrary overconvergent $F$-isocrystals, which then
(by appropriate d\'evissages) should be useful for extending the
stability of overholonomicity under the remaining
cohomological operations. However, these assertions,
as well as Shiho's conjecture, remain in 
the future at the time of this writing.

\section{$p$-adic Hodge theory}

By ``$p$-adic Hodge theory'', we will mean two things: the study of
the interrelated structures on the various cohomologies (notably
de Rham and $p$-adic \'etale) of varieties over $p$-adic fields, and
the study of the abstract objects (Galois representations and various
parameter modules) arising from these
interrelationships. In particular, one central question is to classify
$p$-adic Galois representations which are ``geometrically interesting'',
though the meaning of ``interesting'' has expanded recently
(see Section~\ref{subsec:nongeom}).

For further discussion of much of this material,
we recommend highly Berger's recent survey \cite{berger-dwork}.

\subsection{Geometric $p$-adic Hodge theory}
\label{subsec:geom}

The origin of $p$-adic Hodge theory lies in Grothendieck's proposal
of a ``mysterious functor'' (\textit{foncteur myst\'erieux}) that would
directly relate the $p$-adic \'etale cohomology and the crystalline cohomology
of an algebraic variety over a $p$-adic field. This question in turn
originated in the cohomology in degree 1, where it naturally occurs phrased
in terms of $p$-divisible groups, and was answered by Fontaine
\cite{fontaine-barsotti, fontaine-barsotti2}. 

Fontaine then constructed 
\cite{fontaine-asterisque1, fontaine-asterisque2} a general setup for 
approaching Grothendieck's question, which we now briefly introduce. 
For $K$ a finite extension
of $\QQ_p$, let $G_K$ denote the absolute Galois group of $K$.
By a \emph{$p$-adic representation} of $G_K$, we will always mean
a finite dimensional $\QQ_p$-vector space $V$ equipped with a continuous
$G_K$-action. By a \emph{ring of $p$-adic periods}, we will mean a
topological $\QQ_p$-algebra $\bB$ equipped with a continuous $G_K$-action,
such that $\bB$ is \emph{$G_K$-regular}: if $b \in \bB$ generates
a $G_K$-stable $\QQ_p$-subspace of $\bB$, then $b \in \bB^*$.
In particular, the fixed ring $\bB^{G_K}$ is a field.
For $V$ a $p$-adic representation and $\bB$ a ring of $p$-adic periods,
we define the ``period space''
\[
D_{\bB}(V) = (V \otimes_{\QQ_p} \bB)^{G_K};
\]
it is easily shown to be
a $\bB^{G_K}$-vector space of dimension less than or equal to the
$\QQ_p$-dimension of $V$. If equality holds, we say $V$ is 
\emph{$\bB$-admissible}.

Fontaine exhibited a number of rings
 of $p$-adic periods
\[
\bB_{\crys}, \bB_{\st}, \bB_{\dR}, \bB_{\HT}
\]
the subscripts respectively abbreviating 
``crystalline'', ``semistable'',
``de Rham'', ``Hodge-Tate'';
conversely, one abbreviates ``$\bB_{\crys}$-admis\-sible'' to ``crystalline''
and so on. These conditions get weaker as you move to the right,
so every crystalline representation is semistable, and so on. 
One may insert the condition ``potentially semistable'' after ``semistable'',
for a representation which becomes semistable upon restriction to $G_{K'}$ for
some finite extension $K'$ of $K$.
With this insertion, the reverse
implications all fail with one exception; see Section~\ref{subsec:analytic}.

Let $X$ be a smooth proper variety over $K$; then the $p$-adic \'etale
cohomology $V = H^i_{\et}(X \times_K \overline{K}, \QQ_p)$ is a $p$-adic representation of $G_K$,
and the period rings can be used to extract the ``hidden structure'' on
$H^i_{\et}(X \times_K \overline{K},\QQ_p)$. 
This was first done (not in this framework or terminology) 
by Tate, who considered
the ring $\bB_{\HT} = \oplus_{n \in \ZZ} \CC_p(n)$, where $\CC_p$ is the 
completed algebraic closure of $K$ and $(n)$ denotes a twist by the
$n$-th power of the cyclotomic character. Any Hodge-Tate representation
thus carries a set of numerical invariants: if $V$ is Hodge-Tate, then
$V \otimes_{\QQ_p} \CC_p$ decomposes as a direct sum of copies of various
$\CC_p(n)$. The collection of these integers $n$ is called the set of
\emph{Hodge-Tate weights} of $V$. Tate conjectured that 
$H^i_{\et}(\QQ_p \times_K \overline{K},X)$ is always Hodge-Tate with Hodge-Tate weights in 
$\{0, \dots, i\}$, and that the multiplicity of the weight $j$ is
equal to the Hodge number $h^{j,i-j} = \dim_K H^j(\Omega^i_{X/K}, X)$.
This was established in limited generality by
Bloch and Kato \cite{bloch-kato}, and in full
by Faltings  \cite{faltings-ht}.

To describe the mysterious functor, Fontaine found that Hodge-Tate 
admissibility was not sufficient. Instead, he passed to the class of 
de Rham representations and 
conjectured the following:
\begin{itemize}
\item
the representation $V = H^i_{\et}(X \times_K \overline{K},\QQ_p)$ is always de Rham;
\item
there is a canonical isomorphism of $D_{\dR}(V)$ (which is a vector
space over $\bB_{\dR}^{G_K} = K$) with the de Rham
cohomology of $X$, under which the Hodge filtration is obtained by
a distinguished filtration on $\bB_{\dR}$;
\item
there is a recipe (only specified later; see below) to get 
back from the de Rham cohomology to the \'etale cohomology.
\end{itemize}
This was proved in several stages. First, in case $X$ has good reduction,
Fontaine expected an analogous set of statements involving $\bB_{\crys}$;
these were proved by Fontaine and Messing \cite{fontaine-messing},
Faltings \cite{faltings}, and later again (using $K$-theoretic
techniques) by Nizio\l\ \cite{niziol}.
In case $X$ has semistable reduction, Fontaine and Jannsen
expected a similar set of statements involving $\bB_{\st}$,
specifying the recipe to return from de Rham
to \'etale cohomology;
this was ultimately established by Tsuji \cite{tsuji}, 
drawing on work of numerous authors.\footnote{Mokrane's
MathSciNet review of \cite{tsuji} details amply the swarm of results
that funnel into Tsuji's work.} Finally, to deduce the de Rham
statement, Tsuji actually shows that $V$ is potentially semistable
(i.e., semistable upon restriction to a suitable $G_{K'}$)
by adding to the mix
de Jong's semistable alterations theorem \cite{dejong}.

\subsection{Abstract $p$-adic Hodge theory}
\label{subsec:analytic}

The ``abstract'' aspect of $p$-adic Hodge theory should be thought of
as analogizing the study of abstract Hodge structures, variations of Hodge
structures, and the like, without direct reference to algebro-geometric objects. 
Here an interesting interrelationship emerges
between $p$-adic Galois representations and the $p$-adic differential
equations considered in Section~\ref{subsec:localmono}.

One important result in the abstract theory is the Colmez-Fontaine theorem,
which classifies certain $p$-adic Galois representations in terms of
simple linear algebra data (thus perhaps justifying the use of the term
``Hodge theory'' in the phrase ``$p$-adic Hodge theory'').
To state it, let $K$ be a finite extension of $\QQ_p$ with maximal
unramified subextension $K_0$.
A \emph{$(\phi,N)$-module} is a finite dimensional $K_0$-vector
space $D$ equipped with a Frobenius-semilinear action $\phi: D \to D$,
and a $K_0$-linear map $N: D \to D$ satisfying $N \phi = p \phi N$.
Such an object is \emph{filtered} if it comes with the data of an
exhaustive separated
descending filtration $\Fil^i$ on $L \otimes_{L_0} V$ by $G$-stable
subspaces. (No condition is made concerning compatibility with $\phi$ or $N$.)

If $V$ is a semistable representation of $G_K$, then
$(B_{\st} \otimes_{\QQ_p} V)^{G_K}$ inherits the structure of a
filtered $(\phi,N)$-module by deriving the $\phi$-action, $N$-action, and
filtration from such structures described explicitly on $B_{\st}$.
Moreover, such filtered $(\phi,N)$-modules are distinguished by a certain
numerical property. For a filtered $(\phi,N)$-module $D$ of rank 1, define
the \emph{Newton number} $t_N(D)$ as the valuation of any element $a \in K$
defined by the relation $\phi(d) = ad$ for $d \in D$ nonzero.
Define the \emph{Hodge number} $t_H(D)$ as the largest $i$ such that
$\Fil^i(D) \neq 0$. For a general filtered $(\phi,N)$-module $D$,
define $t_N(D) = t_N(\wedge^{\dim(D)} D)$ and
$t_H(D) = t_H(\wedge^{\dim(D)} D)$.
\begin{theorem} \label{thm:cf}
A filtered $(\phi,N)$-module $D$ arises from a semistable representation
if and only if:
\begin{enumerate}
\item[(a)] $t_N(D) = t_H(D)$;
\item[(b)] for any $(\phi,N)$-submodule $D'$ of $D$ equipped with the induced
filtration, $t_N(D') \geq t_H(D')$.
\end{enumerate}
\end{theorem}
To put this in the formalism of semistability of vector bundles,
if one defines the ``degree'' of $D$ to be $t_H(D) - t_N(D)$, then
$D$ arises from a semistable representation if and only if it is
``semistable of slope $0$''.

Theorem~\ref{thm:cf} was originally proved by Colmez-Fontaine
\cite{colmez-fontaine} by an ingenious but technically difficult\footnote{To
be fair, the same characterization could reasonably be made of
\cite{kedlaya-local}, which underlies Berger's proof; one counterargument
is that the results of \cite{kedlaya-local} are applicable a bit more
broadly.} argument.
A different proof was given by Berger \cite{berger-weak} by working with
$(\phi, \Gamma)$-modules\footnote{A $(\phi,\Gamma)$-module is an algebraic
object that describes a $p$-adic representation by replacing the
complicated Galois action with a simple action on a more complicated ring.
Again, see \cite{berger-dwork} for more on the utility of such objects.}
over the Robba ring, and at one point
invoking Kedlaya's slope filtration theorem (see
Section~\ref{subsec:localmono});
a variant of this argument in the crystalline case was given by Kisin
\cite{kisin-wam}, which as an offshoot established some classification
results for $p$-divisible groups and finite flat group schemes conjectured
by Breuil (unpublished).

Another key result is Fontaine's conjecture that 
every de Rham representation is potentially semistable. (Recall that 
potential semistability for representations arising from \'etale cohomology
was established by Tsuji; see Section~\ref{subsec:geom}.) This 
conjecture was first
proved by Berger \cite{berger-derham} again using 
$(\phi, \Gamma)$-modules over the Robba ring; this time the key input from
that theory is precisely Theorem~\ref{thm:plmt}.
Subsequently, proofs within the ``Fontaine context''
were given by Colmez \cite{colmez-pst} and Fontaine
\cite{fontaine-pst}.

An additional direction, which has some relevance for 
applications\footnote{Example: the upcoming Brandeis thesis of 
Seunghwan Chang will apply these ideas to the formulation of
variants of Serre's
conjecture on modularity of Galois representations.}
in number theory, is the relationship between linear algebraic
descriptions of $p$-adic representations and Galois cohomology.
This was first described by Herr \cite{herr, herr2}, who demonstrated
its utility by recovering some valuable results of Tate, such as local
duality. In a different direction, Marmora
\cite{marmora} has described a relationship between Swan conductors
of a potentially semistable representation and ``irregularity''
of filtered modules; this is an mixed-characteristic analogue of the
Grothendieck-Ogg-Shafarevich formula in rigid cohomology (see
Section~\ref{subsec:overcon}).

\subsection{Nongeometric representations and a Langlands correspondence}
\label{subsec:nongeom}

In recent years, the subclass among $p$-adic Galois representations of those
considered ``geometrically interesting'' has been significantly enlarged.
The old definition would have restricted to those which are potentially
semistable (or de Rham, which is the same; see
Section~\ref{subsec:analytic}), as those are the ones which can occur within
$p$-adic \'etale cohomology. However, the theory of modular forms
suggests that a bigger class should be considered, including certain
``interpolations'' among geometric representations.

Specifically, to address  a qualitative\footnote{The
quantitative version of the Gouv\^{e}a-Mazur conjecture is actually
incorrect as stated \cite{buzzard-calegari}.} version of the
Gouv\^{e}a-Mazur conjecture \cite{gouvea-mazur} on $p$-adic variation
of modular forms, Coleman 
\cite{coleman}
introduced the class of \emph{overconvergent\footnote{This
word will look familiar; it refers to the fact that roughly speaking,
these are defined as sections of a sheaf on a dagger space associated
to a modular curve. The ``roughly'' is because that sheaf may have to be
raised to a non-integral power; Iovita and Stevens are working on several
ways to render Coleman's \emph{ad hoc} workaround for this more systematic.}
 modular forms}.
Coleman and Mazur \cite{coleman-mazur} demonstrated that these were attached
to (global) Galois representations naturally para\-metrized by a rigid
analytic curve (the \emph{eigencurve}).
Kisin \cite{kisin-fm} showed that of these, only the representations
attached to classical modular forms have local representations which
are potentially semistable.

In order to understand the local representations arising on the eigencurve,
e.g., to prove theorems\footnote{This has already begun:
Kisin has given a strengthening \cite{kisin-geom}
of his modularity results using the work of Berger and Breuil.} 
on the modularity of Galois representations in
the vein of Taylor and Wiles \cite{wiles, taylor-wiles}, one would like
to describe a classification of $p$-adic representations in the spirit
of the local Langlands correspondence of Harris and Taylor
\cite{harris-taylor} and Henniart \cite{henniart} (also see
\cite{carayol}) for $\ell$-adic representations. Since there are many
more $p$-adic representations, the $p$-adic correspondence will necessarily
have a different flavor: the ``automorphic'' representations of $\GL_n(K)$,
for $K$ a $p$-adic field, corresponding to $n$-dimensional representations
of $G_K$ should be \emph{infinite-dimensional} vector spaces. The appropriate
such representations appear to be the locally analytic representations
of Schneider and Teitelbaum 
\cite{schneider-teitelbaum1, schneider-teitelbaum2, 
schneider-teitelbaum3, schneider-teitelbaum4}.

Almost all work in this direction so far has focused on the case of
$\GL_2$. Some initial evidence in this case was provided by the work of
Breuil and M\'ezard \cite{breuil-mezard}, and by subsequent work of
Breuil \cite{breuil1, breuil2}. This has led Breuil to a series of
predictions about 2-dimensional representations of $G_K$ and about
integral structures on said representations, which appear to relate to
Banach lattices on the automorphic sides. Some of these predictions in
the case of crystalline representations have now been checked by
Berger and Breuil \cite{berger-breuil1, berger-breuil2}.

As promised above, further investigation has begun to suggest new classes
of representations which are meaningful in the Langlands correspondence.
One such class is the class of \emph{trianguline} representations introduced
by Colmez \cite{colmez-serie}. These are 2-dimensional representations
of $G_K$ whose $(\phi,\Gamma)$-modules over the Robba ring
can be written as extensions of one rank 1 module by another.
(Note that the two rank 1 modules do not in general correspond to
representations; this characterization thus lies firmly within the
theory of Frobenius modules over the Robba ring. See 
Section~\ref{subsec:localmono}.)
Colmez proposes a conjectural correspondence between trianguline
representations and unitary principal series representations of $\GL_2(K)$.
For these representations, Breuil's conjecture on the mod $p$ reduction
has been verified by Berger \cite{berger-dim2}.

\subsection{Towards nonabelian $p$-adic Hodge theory}

We end with a brief mention of some results in the direction of
``nonabelian $p$-adic Hodge theory''. If one thinks of ``abelian''
Hodge theory, $p$-adic or otherwise, as the study of extra
structures on the abelianization of the fundamental group of a 
variety, then nonabelian Hodge theory should be the study of extra
structures on possibly nonabelian quotients of the fundamental group.
For example, in ordinary Hodge theory, Simpson \cite{simpson} gives
an analogue of the Hodge decomposition for the cohomology of a variety
with coefficients in a local system (representation of the fundamental
group), using Higgs bundles. 

One source of inspiration for nonabelian Hodge theory was Deligne's
tract \cite{deligne} on the projective line minus three points; it gives
a natural Tannakian
framework for considering the Hodge theory of \emph{unipotent}
quotients of the fundamental group. In particular, there is a realization
of this framework corresponding to crystalline cohomology, giving rise
to a crystalline fundamental group. 
This generalizes pretty broadly,
to suitable log schemes over a field of positive characteristic
\cite{chiarellotto-lestum, shiho1, shiho2}.

A more robust foundation for nonabelian Hodge theory comes from
rational homotopy theory, as in the work of Katzarkov-Pantev-Toen
\cite{kpt}. This for starters gives a more uniform construction of
crystalline fundamental groups \cite{hain-kim}. That in turn should fit
into a fuller nonabelian $p$-adic Hodge theory parallel to that of
\cite{kpt}; results in this direction have been obtained by
Olsson \cite{olsson-fhomotopy, olsson-nonabelian}.

We cannot conclude without pointing out that something as seemingly
abstruse nonabelian $p$-adic Hodge theory may have concrete applications
to Diophantine equations! Kim \cite{kim-siegel} has suggested a
nonabelian generalization of Chabauty's method \cite{chabauty} for bounding
the number of, and in some cases determining the exact set of,
rational points on a curve over a number field. Chabauty's
method has been rendered practical by a series of refinements
\cite{coleman-chabauty, flynn, flynn-wetherell}, but usually only works
when the Mordell-Weil group of the Jacobian has rank less than the genus
of the curve (so that the group lies within a closed $p$-adic analytic 
subvariety of the Jacobian). 
It is hoped 
that the nonabelian version, which would take place on a higher Albanese
variety \cite{hain-albanese} instead of the Jacobian, 
may yield additional practical
results in cases where this ``Chabauty condition'' is not satisfied.

\bibliographystyle{amsalpha}

\end{document}